\documentclass[11pt]{article}

\usepackage[a4paper,margin=25mm]{geometry}
\usepackage{amsmath,amssymb}
\usepackage{bm}
\usepackage{graphicx}
\usepackage{cite}        
\usepackage{hyperref}    

\hypersetup{
  colorlinks=true,
  linkcolor=blue,
  citecolor=blue,
  urlcolor=blue
}

\title{A Numerical Analysis for Pursuit--Evasion Games\\
under the Stackelberg Equilibrium}

\author{%
  Kazuhiro Horie\\
  Acquisition, Technology and Logistics Agency, Ministry of Defense, Japan%
  \thanks{This manuscript is an English translation of a paper originally published in Japanese in the \textit{Journal of the Japan Society for Aeronautical and Space Sciences}, Vol.~63, No.~5, pp.~204--209 (2015).
At the time of the original publication, the author was affiliated with the Air Systems Research Center, Technical Research and Development Institute, Ministry of Defense, Japan}%
}

\date{} 

\begin{document}
\maketitle


\begin{abstract}
This study presents a numerical analysis framework for pursuit--evasion differential games under the Stackelberg equilibrium.
The Semi-DCNLP method is introduced as an optimization solver for flight path optimization problems formulated under the Stackelberg equilibrium.
In a spacecraft pursuit-evasion differential game, the Semi-DCNLP method produces trajectories that closely match those obtained with a shooting method. This agreement supports the validity of the proposed approach.
\end{abstract}


\vspace{2mm}

\section{Introduction}
Flight path optimization has advanced rapidly since the latter half of the twentieth century, with a close relationship to optimal control theory, aiming to control the trajectories of aircraft and spacecraft so as to minimize a prescribed performance index.

Building on the theoretical foundations of flight path optimization systematized by Bryson and Ho \cite{BrysonHo1975}, numerical methods for solving more practical problems have also been developed. These methods can be broadly classified into two categories: indirect methods, which analytically derive the control variables using the calculus of variations or Pontryagin's minimum principle and transform the problem into a boundary-value problem of ordinary differential equations to be solved numerically; and direct methods, which approximate the control variables by a finite number of parameters and minimize the performance index with respect to these parameters.

In numerical solution methods for flight path optimization, indirect methods were once predominant, since they allow the control variables to be obtained analytically, thereby saving computational resources while achieving high-accuracy solutions. On the other hand, a major drawback of indirect methods is that the resulting boundary-value problems of ordinary differential equations can become difficult to solve numerically when strong nonlinearities are present or when complex boundary conditions involving interior constraints arise. As problems become more realistic, these limitations tend to become increasingly apparent.

Direct methods have been studied primarily by researchers at Boeing since the early 1970s, when indirect methods were at their peak \cite{HargravesEtAl1981}. A major challenge of direct methods lies in how to parameterize the variables used to optimize the performance index. Until the late 1980s, only limited success had been achieved; however, in the late 1980s, a broadly applicable method known as DCNLP (Direct Collocation with Nonlinear Programming) was developed by researchers at Boeing in collaboration with the U.S. Air Force. Although the detailed derivation is referred to in Ref.~\cite{HargravesParis1987}, the DCNLP method transforms the ordinary differential equations describing the equations of motion in a flight path optimization problem into algebraic equations using a numerical technique known as collocation. By treating time-discretized values of the state and control variables as parameters, the original flight path optimization problem is converted into a nonlinear programming problem. As a result, the flight path optimization problem can be reformulated as a constrained nonlinear programming problem consisting of algebraic equations and inequalities.

Although this approach requires substantial computational resources because it minimizes the performance index with respect to a large number of parameters under numerous algebraic equality and inequality constraints, this requirement is no longer a serious drawback in view of recent advances in nonlinear programming solvers and modern computational capabilities. The advantages of the method include the fact that it does not require complicated analytical derivations of control variables based on the calculus of variations or Pontryagin's minimum principle, and that it exhibits good convergence properties even for strongly nonlinear problems. Consequently, the DCNLP method has attracted considerable attention as a means of solving complex problems with relative ease, and successful results have been reported for a wide range of flight path optimization problems \cite{EnrightConway1991,HorieConway2000,YamamotoEtAl2003,IshiiEtAl2008}. Furthermore, at NASA, a trajectory optimization tool known as OTIS \cite{HargravesEtAl1992OTIS}, which applies the DCNLP method, has been developed and put into practical use.

Building on these successes, numerical solution methods based on DCNLP were extended by Horie and Conway \cite{HoriePhD2002,HorieConway2006} to pursuit--evasion differential games involving two players, which represent a class of more complex problems. This approach, referred to as the Semi-DCNLP method, aims to obtain saddle-point equilibrium solutions for pursuit--evasion differential games in which two players determine their strategies simultaneously. In this method, optimality conditions and costate equations for the pursuer (or evader) are added to the original differential game formulation, and the resulting problem is treated as an optimization problem for the evader (or pursuer) to which the DCNLP method is applied. Since its proposal, the Semi-DCNLP method has been successfully applied to obtain numerical solutions for several problems, including interception problems against ballistic missiles \cite{PontaniConway2008,ParishMS2008}. However, because this method does not guarantee that the evader (or pursuer) satisfies the optimality conditions, it cannot be said to ensure the sufficient conditions for a saddle-point equilibrium, and an additional verification of these conditions is required \cite{HoriePhD2002}.

On the other hand, in addition to saddle-point equilibria, differential games also admit solutions under another equilibrium concept known as the Stackelberg equilibrium. Both Stackelberg and saddle-point equilibria assume multiple decision makers with conflicting objectives and no incentive to cooperate. However, Stackelberg equilibria differ in that decision making is sequential. The decision maker who acts later (the follower) determines its strategy after observing the choice made by the earlier decision maker (the leader).

That is, under the assumption that the follower makes an optimal response to the leader’s decision, the leader chooses its strategy so as to optimize its own objective, and the solution of this formulation constitutes a Stackelberg equilibrium. By virtue of this structure, although it is necessary to analytically derive the follower’s optimal response, the solution of a differential game under a Stackelberg equilibrium can ultimately be obtained by solving an optimization problem for the leader that incorporates the follower’s optimal response as a constraint.

The present study addresses pursuit--evasion differential game problems and aims to overcome the aforementioned weakness of the Semi-DCNLP method by introducing the Stackelberg equilibrium concept in place of the saddle-point equilibrium. This approach eliminates the need to separately demonstrate that the obtained solution satisfies the sufficient conditions for a saddle-point equilibrium when using the Semi-DCNLP method, and is therefore expected to facilitate a more straightforward computation of solutions to pursuit--evasion differential games.

\section{Formulation}
In this section, the pursuit--evasion differential game problem under a Stackelberg equilibrium, which is the focus of the present study, is formulated. First, the equations of motion of the pursuer and the evader, together with their initial conditions, are defined as follows:
\begin{align}
\frac{d\bm{x_p}}{dt} &= \bm{f_p}(\bm{x_p},\bm{u_p},t) \tag{1}\\
\frac{d\bm{x_e}}{dt} &= \bm{f_e}(\bm{x_e},\bm{u_e},t) \tag{2}\\
\bm{x_p}(t_0) &= \bm{x_{p0}} \tag{3}\\
\bm{x_e}(t_0) &= \bm{x_{e0}} \tag{4}
\end{align}
For simplicity, the equations of motion of the pursuer and the evader are assumed to depend only on their respective state variables $\bm{x_p}, \bm{x_e}$ and control variables $\bm{u_p}, \bm{u_e}$. This assumption is generally valid when a fixed reference frame is adopted.

Next, a terminal condition is specified:
\begin{align}
\bm{\psi}(\bm{x_p}(t_f), \bm{x_e}(t_f),t_f)=0 \tag{5}
\end{align}
The game terminates when the terminal condition is satisfied, and the corresponding time is defined as the terminal time $t_f$. In the present study, no constraints on the state variables or control variables are considered.

A zero-sum game is assumed in this study, and the corresponding performance index (cost functional) $J$ is defined as
\begin{align}
J = \bm{\phi}(\bm{x_p}(t_f), \bm{x_e}(t_f),t_f) + \int_{t_0}^{t_f} L(\bm{x_p},\bm{x_e},\bm{u_p},\bm{u_e},t)\,dt. \tag{6}
\end{align}

In differential games, feedback strategies are introduced to determine the control variables as functions of the state variables. In this study, we use an open-loop representation of feedback strategies, i.e., a representation that specifies the optimal trajectories for the pursuer and the evader as functions of their initial state variables, and express it as
\begin{align}
\bm{u_p} &= \bm{\gamma_p}(\bm{x_{p0}},\bm{x_{e0}},t) \tag{7}\\
\bm{u_e} &= \bm{\gamma_e}(\bm{x_{p0}},\bm{x_{e0}},t). \tag{8}
\end{align}

As described in the Introduction, the Stackelberg equilibrium considered in this study divides the players into a leader and a follower, where the leader selects its optimal strategy under the assumption that the follower will choose an optimal strategy, namely, optimal control variables.

Either the pursuer or the evader is assumed to be the leader, with the other acting as the follower. By denoting the leader and follower with subscripts $l$ and $f$, respectively, the value of the game $V$ is expressed as follows. The leader is assumed to choose its strategy so as to maximize the performance index, whereas the follower seeks to minimize it:
\begin{align}
V=\max_{\bm{u_l}}\min_{\bm{u_f}} J. \tag{9}
\end{align}

A system of equations for differential games under a Stackelberg equilibrium was derived by Basar and Olsder \cite{BasarOlsder1999} as an extension of the corresponding equations for standard optimal control problems. Applying this system to Eqs.~(1)--(9) yields the following results. First, by fixing the leader’s control variables and determining the follower’s optimal trajectory, the following equations must be satisfied.

\noindent\textbf{Hamiltonian:}
\begin{align}
H = L(\bm{x_l},\bm{x_f},\bm{u_l},\bm{u_f},t)+ \bm{\lambda_l}^{T} \bm{f_l}(\bm{x_l},\bm{u_l},t)+ \bm{\lambda_f}^{T} \bm{f_f}(\bm{x_f},\bm{u_f},t). \tag{10}
\end{align}

\noindent\textbf{Follower optimality condition:}
\begin{align}
\bm{u_f} = \arg\min_{\bm{u_f}} H. \tag{11}
\end{align}

\noindent\textbf{Costate equations:}
\begin{align}
\frac{d\bm{\lambda_l}}{dt} &= -\left(\frac{\partial H}{\partial \bm{x_l}}\right)^{T} \tag{12}\\
\frac{d\bm{\lambda_f}}{dt} &= -\left(\frac{\partial H}{\partial \bm {x_f}}\right)^{T}. \tag{13}
\end{align}

\noindent\textbf{Terminal conditions for costates:}
\begin{align}
\bm{\lambda_l}(t_f) &= \left(\frac{\partial \phi}{\partial \bm{x_l}}\right)^{T}+\left(\frac{\partial \bm{\psi}}{\partial \bm{x_l}}\right)^{T}\bm{\nu} \tag{14}\\
\bm{\lambda_f}(t_f) &= \left(\frac{\partial \phi}{\partial \bm{x_f}}\right)^{T}+\left(\frac{\partial \bm{\psi}}{\partial \bm{x_f}}\right)^{T}\bm{\nu} \tag{15}
\end{align}

\noindent\textbf{Transversality condition:}
\begin{align}
\left[H+\frac{\partial \phi}{\partial t}+\bm{\nu}^{T}\left(\frac{\partial \bm{\psi}}{\partial t}\right)\right]_{t=t_f}=0. \tag{16}
\end{align}

The state variables of the follower and the leader that satisfy Eqs.~(1)--(5) and (11)--(16) represent the respective state variables under the follower's optimal strategy. Note that the leader’s costate variable $\lambda_l$ does not appear at times other than the terminal time except through Eq.~(12); therefore, a solution can be obtained without explicitly solving the leader’s costate equation in Eq.~(12). Accordingly, it suffices to require that the state variables of the follower and the leader satisfy Eqs.~(1)--(5), (11), and (13)--(16); if the leader’s control variables are then determined so as to optimize the performance index $J$ under these conditions, the Stackelberg equilibrium solution of the given problem is obtained.

That is, the pursuit--evasion differential game problem under a Stackelberg equilibrium is reduced to solving the following problem subject to the constraints given by Eqs.~(1)--(5), (11), and (13)--(16):
\begin{align}
V=\max_{\bm{u_l}}J=\min_{\bm{u_l}}(-J). \tag{17}
\end{align}

\section{Application of the Semi-DCNLP Method to Problems under a Stackelberg Equilibrium}
In the DCNLP method, which is based on nonlinear programming, it is not possible to simultaneously maximize the performance index with respect to one player while minimizing it with respect to the other. Consequently, the method cannot be directly applied to differential game problems under a saddle-point equilibrium.

Therefore, as described in Section~1, the Semi-DCNLP method enables the use of nonlinear programming by augmenting the original pursuit--evasion differential game problem with the optimality conditions and costate equations for either the pursuer or the evader \cite{HoriePhD2002,HorieConway2006}. However, according to Horie \cite{HoriePhD2002}, in order to ensure that the solution obtained by this approach indeed corresponds to a saddle-point equilibrium, the additional costate variables must coincide with the Lagrange multipliers obtained in the nonlinear programming solution process. This requirement necessitates a post-processing procedure based on the method proposed by Enright and Conway \cite{EnrightConway1992}.

On the other hand, as formulated in Section~2, when the equilibrium concept of the differential game is changed from a saddle-point equilibrium to a Stackelberg equilibrium, the follower’s optimal response is fully characterized by the differential equations, algebraic equations, and algebraic inequalities given in Eqs.~(11) and (13)--(16). Consequently, by discretizing Eqs.~(1)--(5), (11), and (13)--(16) at appropriate time intervals and replacing the differential equations with algebraic equations, the problem can ultimately be reduced---through a process analogous to that of the DCNLP method---to a constrained nonlinear programming problem that optimizes the leader’s performance index.

Although this procedure is identical to that used in the Semi-DCNLP method for solving differential games under a saddle-point equilibrium, the resulting numerical solution satisfies the Stackelberg equilibrium conditions. As a result, the additional step required in the saddle-point case---namely, verifying that the costate variables coincide with the Lagrange multipliers obtained from the nonlinear program---can be avoided.

\section{Numerical Example}
In this section, the procedure described up to Section~3 is applied to a simple problem. The validity of the solution obtained by the Semi-DCNLP method is examined by comparing it with a solution obtained by an indirect method, in which the problem is formulated as a two-point boundary-value problem and solved numerically using a shooting method.

Well-known pursuit--evasion differential games include the Homicidal Chauffeur Problem and the Dorichobrachistochrone Problem, both introduced by Isaacs \cite{Isaacs1965}. On the other hand, solutions for these problems have been obtained under saddle-point equilibria, and solutions under a Stackelberg equilibrium are not necessarily available. Moreover, because the control variables in these problems are constrained, it is somewhat complicated to treat them as two-point boundary-value problems via an indirect method and to obtain numerical solutions using a shooting method.

Therefore, a spacecraft pursuit--evasion problem is adopted as a benchmark problem in this study. The spacecraft pursuit--evasion problem is also a classical problem that was addressed relatively early---in the mid-1970s---by Anderson and Grazier \cite{GrazierMS1973,AndersonGrazier1976}. Since the control variables are unconstrained, solutions can be obtained relatively easily even by a shooting method; for this reason, it is selected as the benchmark problem.

\begin{figure}[ht]
  \centering
  \includegraphics[width=0.78\linewidth]{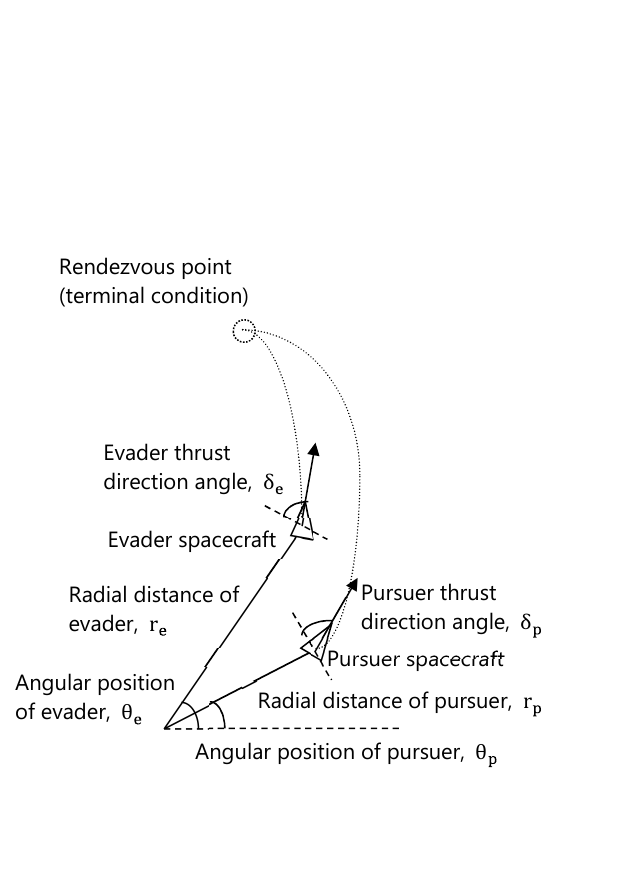}
  \caption{Schematic overview of the spacecraft pursuit--evasion problem.}
  \label{fig:fig1}
\end{figure}

As shown in Fig.~\ref{fig:fig1}, this problem considers two point masses moving on a two-dimensional plane under a gravitational field. The game terminates when the pursuer captures the evader. The pursuer controls its thrust direction angle to minimize the capture time, whereas the evader controls its thrust direction angle to maximize the time until capture. The pursuer’s thrust is set larger than the evader’s thrust so that capture occurs and the game admits a well-defined solution.

The gravitational field is assumed to be inversely proportional to the square of distance, and atmospheric effects are neglected, yielding a simplified yet reasonably realistic model. In addition, the spacecraft mass is assumed to remain constant, neglecting propellant consumption. Furthermore, the two spacecraft are assumed to have the same mass, so that the equations of motion can be normalized by the spacecraft mass. The resulting equations of motion are given as follows:
\begin{align}
\frac{dv_{rp}}{dt} &= T_p\sin\delta_p-\frac{\mu}{r_p^{2}}+\frac{v_{\theta p}^{2}}{r_p} \tag{18}\\
\frac{dv_{\theta p}}{dt} &= T_p\cos\delta_p-\frac{v_{rp}v_{\theta p}}{r_p} \tag{19}\\
\frac{dr_p}{dt} &= v_{rp} \tag{20}\\
\frac{d\theta_p}{dt} &= \frac{v_{\theta p}}{r_p} \tag{21}\\
\frac{dv_{re}}{dt} &= T_e\sin\delta_e-\frac{\mu}{r_e^{2}}+\frac{v_{\theta e}^{2}}{r_e} \tag{22}\\
\frac{dv_{\theta e}}{dt} &= T_e\cos\delta_e-\frac{v_{re}v_{\theta e}}{r_e} \tag{23}\\
\frac{dr_e}{dt} &= v_{re} \tag{24}\\
\frac{d\theta_e}{dt} &= \frac{v_{\theta e}}{r_e}. \tag{25}
\end{align}

The game terminates when the following conditions are satisfied:
\begin{align}
r_p=r_e \tag{26}\\
\theta_p=\theta_e. \tag{27}
\end{align}

As stated above, the performance index is taken as the time required for the pursuer to capture the evader. When the evader is the leader and the pursuer is the follower in the Stackelberg equilibrium, the value of the game is written as
\begin{align}
V=\max_{\delta_e}\min_{\delta_p} t_f. \tag{28}
\end{align}

Using Eqs.~(10), (11), and (13)--(16), the conditions required to obtain the Stackelberg equilibrium solution are derived as follows.

\noindent\textbf{Hamiltonian:}
\begin{align}
\begin{split}
H =\ &1
+\lambda_{v_{rp}}\left(T_p\sin\delta_p-\frac{\mu}{r_p^{2}}+\frac{v_{\theta p}^{2}}{r_p}\right)
+\lambda_{v_{\theta p}}\left(T_p\cos\delta_p-\frac{v_{rp}v_{\theta p}}{r_p}\right)\\
&+\lambda_{r_p}v_{rp}+\lambda_{\theta_p}\frac{v_{\theta p}}{r_p}
+\lambda_{v_{re}}\left(T_e\sin\delta_e-\frac{\mu}{r_e^{2}}+\frac{v_{\theta e}^{2}}{r_e}\right)\\
&+\lambda_{v_{\theta e}}\left(T_e\cos\delta_e-\frac{v_{re}v_{\theta e}}{r_e}\right)
+\lambda_{r_e}v_{re}+\lambda_{\theta_e}\frac{v_{\theta e}}{r_e}.
\end{split}
\tag{29}
\end{align}

\noindent\textbf{Pursuer optimality condition:}

Since there is no constraint on the pursuer’s control variable, namely the thrust direction angle $\delta_p$, the optimality condition (11) yields
\begin{align}
\frac{\partial H}{\partial \delta_p}
&=\lambda_{v_{rp}}T_p\cos\delta_p-\lambda_{v_{\theta p}}T_p\sin\delta_p=0 \tag{30}\\
\frac{\partial^{2}H}{\partial \delta_p^{2}}
&=-\lambda_{v_{rp}}T_p\sin\delta_p-\lambda_{v_{\theta p}}T_p\cos\delta_p>0. \tag{31}
\end{align}

\noindent\textbf{Pursuer costate equations:}
\begin{align}
\frac{d\lambda_{v_{rp}}}{dt} &= \lambda_{v_{\theta p}}\frac{v_{\theta p}}{r_p}-\lambda_{r_p} \tag{32}\\
\frac{d\lambda_{v_{\theta p}}}{dt} &= \frac{-2\lambda_{v_{rp}}v_{\theta p}+\lambda_{v_{\theta p}}v_{rp}-\lambda_{\theta_p}}{r_p} \tag{33}\\
\frac{d\lambda_{r_p}}{dt} &= \frac{-2\lambda_{v_{rp}}\mu+\lambda_{v_{rp}}v_{\theta p}^{2}r_p-\lambda_{v_{\theta p}}v_{rp}v_{\theta p}r_p+\lambda_{\theta_p}v_{\theta p}r_p}{r_p^{3}} \tag{34}\\
\frac{d\lambda_{\theta_p}}{dt} &= 0 \tag{35}
\end{align}

\noindent\textbf{Terminal conditions for the costates:}
\begin{align}
\lambda_{v_{rp}}(t_f)=0 \tag{36}\\
\lambda_{v_{\theta p}}(t_f)=0 \tag{37}\\
\lambda_{r_p}(t_f)=\nu_{1} \tag{38}\\
\lambda_{\theta_p}(t_f)=\nu_{2} \tag{39}\\
\lambda_{v_{re}}(t_f)=0 \tag{40}\\
\lambda_{v_{\theta e}}(t_f)=0 \tag{41}\\
\lambda_{r_e}(t_f)=-\nu_{1} \tag{42}\\
\lambda_{\theta_e}(t_f)=-\nu_{2} \tag{43}
\end{align}

\noindent\textbf{Transversality condition:}
\begin{align}
\left[
1+\lambda_{r_p}(v_{rp}-v_{re})
+\lambda_{\theta_p}\left(\frac{v_{\theta p}}{r_p}-\frac{v_{\theta e}}{r_e}\right)
\right]_{t=t_f}=0. \tag{44}
\end{align}

Accordingly, the following problem is solved subject to the constraints in Eqs.~(18)--(27) and (30)--(44):
\begin{align}
V=\max_{\delta_e}t_f=\min_{\delta_e}(-t_f). \tag{45}
\end{align}

The example problem is solved using the Semi-DCNLP method with the following initial conditions:
\begin{align}
v_{rp}=0,\ v_{\theta p}=0,\ r_p=1.0,\ \theta_p=0,\ v_{re}=0,\ 
 v_{\theta e}=0.9759,\ r_e=1.05,\ \theta_e=0.4. \tag{46}
\end{align}
The pursuer thrust $T_p$, the evader thrust $T_e$, and the gravitational constant $\mu$ are set as follows:
\begin{align}
T_p=0.05,\ T_e=0.0025,\ \mu=1.0. \tag{47}
\end{align}

Although the Hermite rule \cite{EnrightConway1992} is commonly used to select collocation points in the Semi-DCNLP method, in this study the Gauss--Lobatto rule proposed by Harman and Conway \cite{HarmanConway1996} is employed to obtain higher accuracy. As the nonlinear programming solver, SNOPT implemented as E04UGF in the NAG FORTRAN Library is used.

\begin{figure}[ht]
  \centering
  \includegraphics[width=0.78\linewidth]{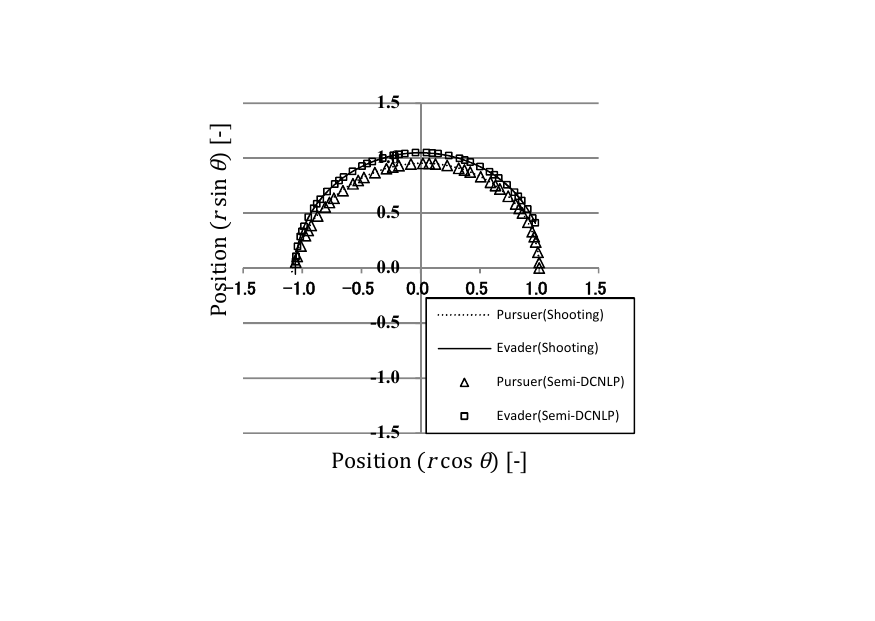}
  \caption{Pursuit--evasion trajectories on the two-dimensional plane.}
  \label{fig:fig2}
\end{figure}

\begin{figure}[ht]
  \centering
  \includegraphics[width=0.78\linewidth]{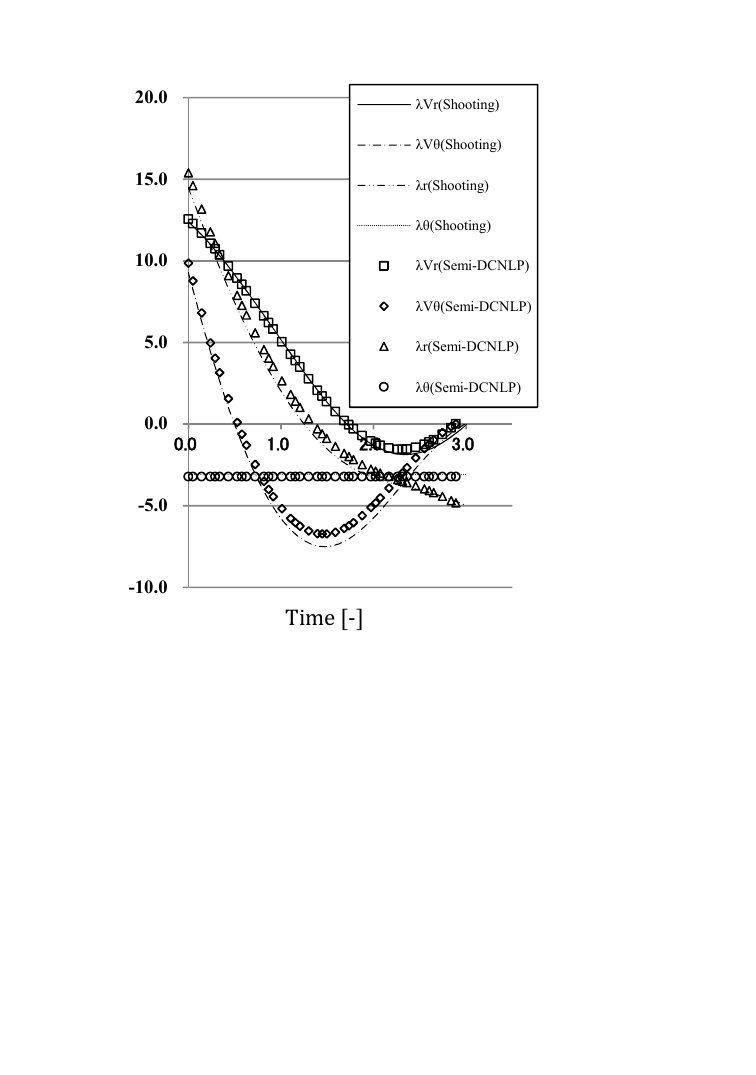}
  \caption{Time histories of the pursuer’s costate variables.}
  \label{fig:fig3}
\end{figure}

\begin{figure}[ht]
  \centering
  \includegraphics[width=0.78\linewidth]{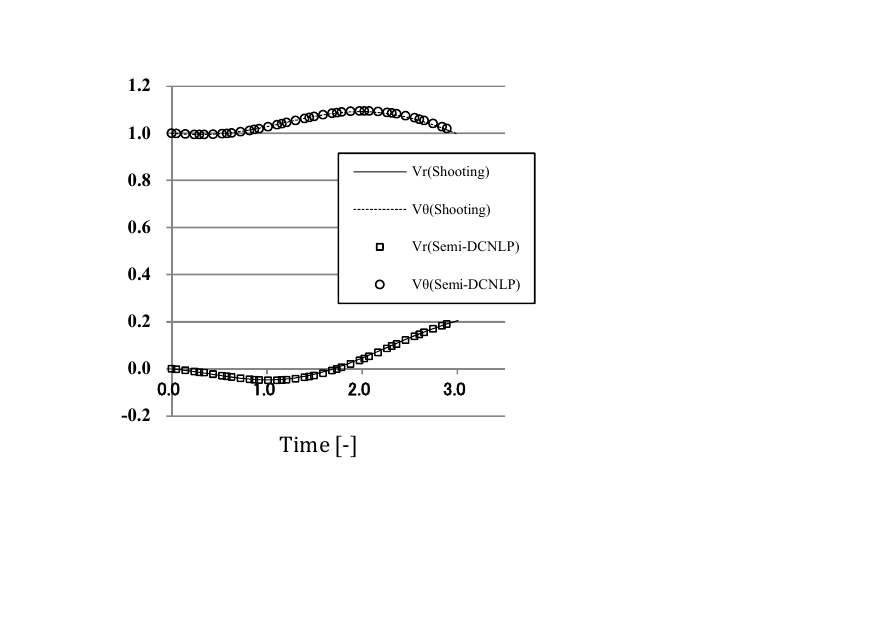}
  \caption{Time histories of the pursuer’s velocity.}
  \label{fig:fig4}
\end{figure}

\begin{figure}[ht]
  \centering
  \includegraphics[width=0.78\linewidth]{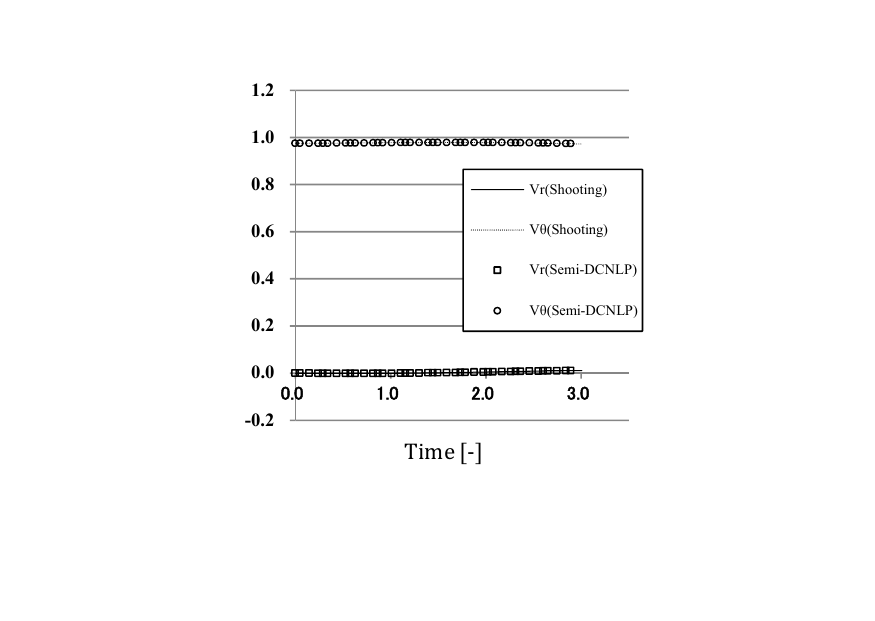}
  \caption{Time histories of the evader’s velocity.}
  \label{fig:fig5}
\end{figure}

\begin{figure}[ht]
  \centering
  \includegraphics[width=0.78\linewidth]{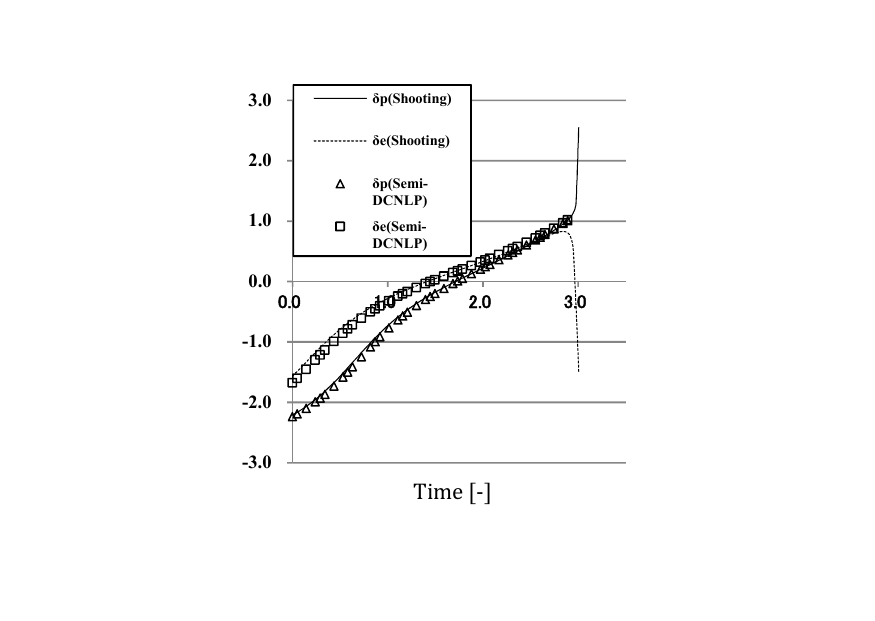}
  \caption{Time histories of the control variables for the pursuer and the evader.}
  \label{fig:fig6}
\end{figure}

The results obtained by the Semi-DCNLP method are shown in Figs.~\ref{fig:fig2}--\ref{fig:fig6}. For validation of the Semi-DCNLP method, these figures also include the results obtained for the same example using the shooting method. In the figures, the Semi-DCNLP results are shown as discrete values at the collocation points, whereas the shooting-method results are shown as continuous curves.

Since the obtained trajectories correspond to the Stackelberg equilibrium solution with the evader as the leader, the evader’s trajectory would be expected to maximize the time until capture by the pursuer. As shown in Fig.~\ref{fig:fig6}, the evader’s time history exhibits a similar trend to that of the pursuer, and the evader’s thrust angle becomes identical to the pursuer’s thrust angle immediately before capture. These observations indicate that the evader applies control so as to prolong the capture time as much as possible.

Near the terminal time, the thrust angles obtained by the shooting method show sharp changes. These changes are considered to arise because the optimality conditions become singular at the terminal point.

In Fig.~\ref{fig:fig2}, both the pursuer and the evader initially move into the interior of the initial orbit and then gradually increase the orbital radius. Capture eventually occurs outside the initial orbit. Similar trajectories have also been obtained in minimum-time rendezvous trajectories of continuous low-thrust spacecraft (e.g., electric propulsion) with near-Earth asteroids \cite{Conway1999AAS}, and can be regarded as a type of trajectory that may arise in minimum-time orbital transfers.

Furthermore, for both the trajectories in Fig.~\ref{fig:fig2} and the time histories of various variables in Figs.~\ref{fig:fig3}--\ref{fig:fig6}, it can be seen that the results of the Semi-DCNLP method are almost identical to those of the shooting method. The terminal time, which corresponds to the game value, is 2.89 for the Semi-DCNLP method, whereas it is 3.01 for the shooting method, resulting in an error of approximately 4\%.

From the above observations, it follows that the qualitative characteristics of the obtained trajectories are consistent with the evader’s objective of maximizing the time until capture in the Stackelberg formulation, and that the solution obtained by the Semi-DCNLP method closely matches the solution obtained by the shooting method. Therefore, the Semi-DCNLP method is shown to be an appropriate approach for numerical analysis of differential game solutions under a Stackelberg equilibrium.

\section{Conclusions}
In this study, it has been demonstrated that solutions to pursuit--evasion differential games under a Stackelberg equilibrium can be obtained by adapting the DCNLP method, which is a direct numerical approach.

Moreover, while numerical solution methods for differential games had traditionally relied primarily on indirect approaches, Horie and Conway developed a method known as the Semi-DCNLP method, which partially incorporates direct methods and enables numerical solutions to be obtained more easily than before. The present study shows that, by adopting the Stackelberg equilibrium instead of the saddle-point equilibrium as the solution concept, the limitations associated with the use of the Semi-DCNLP method can be further relaxed.

It should be noted that, in the present study, no constraints are imposed on the control variables or the state variables. Although this assumption significantly restricts applicability to realistic problems, constraints on the leader’s control variables and state variables, as well as on the follower’s state variables, can be expected to be readily incorporated by formulating the problem as a constrained nonlinear programming problem \cite{HargravesParis1987}. On the other hand, the follower’s control variables are obtained analytically from the optimality conditions, and therefore cannot be treated in the same manner. Nevertheless, it is considered possible to address this issue---albeit with some loss of accuracy---by introducing slack variables at the formulation stage and replacing constrained control variables with new unconstrained control variables.

In the future, numerical analyses of complex and realistic differential game problems are expected to become feasible in various fields by adopting the Stackelberg equilibrium as the solution concept and applying the Semi-DCNLP method, which exhibits favorable convergence properties. This may broaden the range of problems that have been difficult to analyze using conventional approaches.


\end{document}